\documentclass[twocolumn, letterpaper, 10 pt, conference]{ieeeconf}

\IEEEoverridecommandlockouts                              
\usepackage[english]{babel}

\IEEEoverridecommandlockouts                              
\usepackage[english]{babel}
\IEEEoverridecommandlockouts                              
\usepackage[english]{babel}
\usepackage{cite}
\usepackage[cmex10]{amsmath}
\usepackage{dsfont}
\usepackage[tight,footnotesize]{subfigure}

\newif\iflibroproyecto\libroproyectofalse
\input{MacrosMatematicas}
\usepackage[mathscr]{eucal}%
\usepackage{tikz}
\usepackage{multirow}
\usepackage{amsfonts}
\usepackage{amssymb,color}
\usepackage{epstopdf}

\usepackage[absolute]{textpos}
\usepackage{bm}

\usepackage{pstricks}
\usepackage{graphicx}
\usepackage{balance}
\usepackage{tikz}
\usetikzlibrary{shapes,arrows}

\usepackage[prependcaption,colorinlistoftodos]{todonotes}

\begin{document}
\title{\LARGE \bf Stability Via Adversarial Training of Neural Network Stochastic Control of Mean-Field Type}


\author{Julian Barreiro-Gomez~~~~~~Salah Eddine Choutri~~~~~~Boualem Djehiche
\thanks{Julian Barreiro-Gomez and Salah Eddine Choutri are with NYUAD Research Institute, New York University Abu Dhabi, PO Box 129188, Abu Dhabi, United Arab Emirates. (e-mails: {\tt\scriptsize jbarreiro@nyu.edu, choutri@nyu.edu).}} 
\thanks{Boualem Djehiche works at the Department of Mathematics, KTH, Stockholm, Sweden. (e-mails: {\tt\scriptsize boualem@kth.se).} }
\thanks{We gratefully acknowledge support from Tamkeen under the NYU Abu Dhabi Research Institute grant CG002. We thank Hatem Hajri for the discussions on adversarial training.
}
}
\maketitle

\begin{abstract}
In this paper, we present an approach to neural network mean-field-type control and its stochastic stability analysis by means of adversarial inputs (aka adversarial attacks). This is a class of data-driven mean-field-type control where the distribution of the variables such as the system states and control inputs are incorporated into the problem. Besides, we present a methodology to validate the feasibility of the approximations of the solutions via neural networks and evaluate their stability. Moreover, we enhance the stability by enlarging the training set with adversarial inputs to obtain a more robust neural network.
Finally, a worked-out example based on the linear-quadratic mean-field type control problem (LQ-MTC) is presented to illustrate our methodology.
\end{abstract}

\begin{keywords}
Neural networks, data-driven control, stability, robustness, supervised machine learning, adversarial training
\end{keywords}

\vspace{-0.3cm}
\section{Introduction}


 Mean-field type control is a topic that attracted a lot of attention since the introduction of mean-field games by Lasry and Lions in their seminal work \cite{Lions} and by Caines, Huang and Malhamé in \cite{Huang}. Andersson and Djehiche in \cite{Djehicheb} introduced a stochastic mean-field type control problem in which the state dynamics and the performance criterion depend on the moments of the state, see also \cite{buckdahn2011general} and \cite{li2012stochastic}. Carmona and Delarue \cite{carmona}, and Buckdahn \textit{et al.}\cite{buckdahn2016stochastic}, later, generalized the problem to include the probability law of the state dynamics. 
 For applications related to mean-field type control and games problems we cite, among many others, \cite[Chapter 16]{CRC_Book_BaTe}, \cite{Games2020,price2022, IFAC_Boualem} and the references therein. 

This class of problems is non-conventional since both the evolution of the state and often the performance functional are influenced by terms that are not directly related to the state or to the control of the decision maker. In a sense, they model a very large number of agents behaving, all, similarly to a representative agent. The latter is impacted by the aggregation of all agents due to the large number. The aggregation effect is modelled as a mean-field term such as, among others, the law of the state, the expectation of the state or its variance.

Solving this problem, analytically is rather challenging as there are no general analytic methods for this purpose.
Therefore the use of numerical methods is often needed to provide approximations to the solutions and several methods have been suggested for finite horizon mean-field type control problems (see e.g. \cite{achdou2016mean,pfeiffer2017numerical}). Furthermore, the recent progress on machine learning technologies made it easier to test and provide more efficient approximations of the solutions to complex mean-field type control problems.

The link between deep learning and mean-field type control and games was recently studied by, among others, Lauri{\`e}re, Carmona and Fouque in series of papers (see e.g. \cite{Fouque2020, carmona2019convergence, carmonaandmathieu, Laurire2021}), where the authors (jointly and/or independently) proposed algorithms for the solution of mean-field type optimal control problems based on approximations of the theoretical solutions by neural networks, using the software package
TensorFlow with its `Stochastic Gradient Descent' optimizer designed for machine learning. However, the stability of neural networks associated to the mean-field type control problems was not considered in the literature so far.

In deep learning, there is an increasing interest in studying and improving the robustness and stability of the trained neural networks see e.g. \cite{Attacks_NN}, where it has been reported that a simple modification in the input data might fool a well-trained neural network, returning a wrong output. For instance, a picture that is previously well-classified by a trained neural network could be incorrectly classified once we perturb one or more pixels in it.
Such perturbations are known as adversarial attacks and can help to characterize how robust and stable a network is. The contribution of this paper is summarized in three points: training, stability evaluation, and stability improvement. 

\textit{Training:} we present an indirect and simple method to train neural networks to learn optimal controls based on data. Inspired by the work in \cite{carmona2019convergence} and \cite{carmonaandmathieu}, we first illustrate how data can be generated by computationally solving a finite-time horizon optimal control problem with decision variables given by the output of the neural network. Then, we design a data-driven (model-free) mean-field-type control using neural networks in a supervised learning fashion. The idea is to design an ofﬂine controller (once trained, only a simple forward run is required) that is more time-efficient than the conventional online optimization-based control approach, which can be time consuming depending on the complexity of the problem. In real life, one can use the data provided by a traffic application such as google maps to train the neural network to give optimal paths.


\textit{Stability:} we borrow the idea of an adversarial attack from the topic of image classification and draw an analogy in the context of stochastic dynamical systems. An adversarial attack, in our sense, is an initial condition that might make the closed-loop neural network control system unstable. 
This concept enabled us to study the stability of the neural network mean-field-type control and empirically characterize the corresponding forward invariant basin of attraction for the closed-loop system composed of the stochastic dynamics and the optimal control/strategies. 

\textit{Stability improvement:}
we improve the stability of the neural network mean-field-type control by enlarging the training set using adversarial data (attacks) generated from the previous phase (Stability). We compare and discuss the resulting data-enhanced closed-loop system and a suitably modified neural network architecture which potentially enhances the stability of the closed-loop system.


The reminder of this paper is organized as follows. In section II, we formulate the mean-field type control problem. In section III, we present an approach that solves first an optimization problem taking as decision variables the output of the neural network, and then it solves the neural network training. In section IV, we define the stochastic stability concept for our mean-field type neural network by means of adversarial inputs or attacks to the closed-loop neural network. We illustrate our stability results in Section V through numerical examples. Section VI concludes the paper.

\section{Mean-Field-Type Control Problem}
We consider a finite-horizon stochastic control problem where the state process is governed by a
stochastic differential equation (SDE) of mean-field type. The drift here 
depends on the state and control as well as their respective probability laws.
For a fixed time horizon $T>0$, let $(\mathbf{\Omega},\mathcal{F},(\mathcal{F}_t)_{0\le t\le T},\mathbf{P})$ be a filtered probability space satisfying the usual conditions, on which we define a standard Brownian motion $B:=(B(t),\  0\le t\le T)$. We assume that $\mathbf{F}:=(\mathcal{F}_t,\ 0\le t\le T)$ is the natural filtration of $B$ augmented by $\mathbf{P}$-null
sets of $\mathcal{F}$. The action space, $U$, is a non-empty, closed and convex subset of $\mathbb{R}$, and $\mathcal{U}$ is the class of measurable, $\mathbf{F}$-adapted and square integrable processes taking values in $U$. 
For any control $u \in \mathcal{U}$, we consider the following SDE
\begin{align}
\label{eq:sde}
\left\{ 
\begin{array}{l}
\displaystyle
dx(t)=f(t,x(t),\mathbb{E}[x(t)],u(t), \mathbb{E}[u(t)]) dt \\
\displaystyle~~~~~~~ + \sigma dB(t), \\ 
\displaystyle
x(0) = x_0,~ x_0\sim 
\mu_0,
\end{array}
\right.
\end{align}
where,
$f: [0,T ] \times \mathbb{R} \times \mathbb{R} \times U \times \mathbb{R} \rightarrow \mathbb{R}$, $\sigma >0.$
The expected cost is given by 
\begin{align} \label{cost}
\Tilde{J}(u) &= \mathbb{E} [\int_{0}^{T} \ell(t,x(t),\mathbb{E}[x(t)],u(t),\mathbb{E}[u(t)]) dt \\ & + \psi(x(T),\mathbb{E}[x(T)]) ], \notag
\end{align}
where,
$\ell : [0,T ] \times \mathbb{R} \times \mathbb{R} \times U \times \mathbb{R} \rightarrow \mathbb{R},$
 and $\psi : \mathbb{R} \times \mathbb{R} \rightarrow \mathbb{R}.$
The mean-field-type control problem is as follows:
\begin{align*}
\mathcal{P}_{\mathrm{MFTC}}:=
\left\{
\begin{array}{l}
\displaystyle \min_{u \in \mathcal{U}}~ \Tilde{J}(u),\\
\text{subject to}~\eqref{eq:sde},~\text{and}~x_0\sim 
\mu_0. 
\end{array}
\right.
\end{align*}
Next, we present the proposed approach to solve  $\mathcal{P}_{\mathrm{MFTC}}.$



%
%
%
%
%
%

\section{Neural Networks for Mean-Field-Type Problems}
\label{sec:NN}

In this section we define, rigorously, what we mean by a  neural network, then we show how it can be solved our mean-field control problem. 
A neural network is usually defined by an architecture, which is essentially, the number of hidden layers, the number of neurons per layer and the activation functions. We define the set of layer functions with input dimension $d_1$, output dimension $d_2$, and activation functions $h_j:\mathbb{R} \rightarrow \mathbb{R}$, $j \in \{0, \dots,n\},$ by
\begin{align*}
    \mathbb{L}^{h_j}_{d_1,d_2}=\{&\phi: \mathbb{R}^{d_1} \rightarrow \mathbb{R}^{d_2} \vert \exists b \in \mathbb{R}^{d_2}, \exists W \in \mathbb{R}^{d_2 \times d_1}, \\ &\forall i \in \{1, \dots, d_2  \}, \phi(x)^i = h_j(b_i+ \underset{k=1}{\sum} W_{ik} z_k)  \},
\end{align*} 
and we denote the set of neural networks with $n$ hidden layers and one output layer by
\begin{align*}
    \mathcal{U}_{\mathbf{N}\mathbf{N}} = \{&g_{\theta}: \mathbb{R}^{d_0} \rightarrow \mathbb{R}^{d_{n+1}} \vert \forall j \in \{0\dots, n \}, \exists \phi^j \in \mathbb{L}^{h_j}_{d_j,d_{j+1}}, \\ &g_{\theta}=\phi^{n} \circ \phi^{n-1}  \circ \dots \circ  \phi^0 \}. 
\end{align*} 
The vector $b$ and matrix $W$ for each layer, are called the parameters of the neural network, which we usually seek to optimize through training. We denote them by $$\theta := \{ W^{(0)},b^{(0)},W^{(1)},b^{(1)},\dots, W^{(n-1)},b^{(n-1)} \}, \forall n \in \mathbb{N}_{>2}, $$ and we denote their set by $\Theta.$

\subsection{Neural Network Training as Data-Driven Control}

We first discretize time as follows. For a finite $T > 0$ and $N_T \in \mathbb{N}^+$, let $\Delta t = T/N_T$ and $t_k = k \Delta t,$ \ $k \in \{0,\dots,N_T-1 \}$. 
%
%
The discretized version of problem $\mathcal{P}_{\mathrm{MFTC}}$ is given by the state dynamics
\begin{align}
x(t&_{k+1}) = x(t_{k}) \label{eq:sde_dis_control}\\
&+ f(x(t_{k}),\mathbb{E}[x(t_{k})],u(t_{k}),\mathbb{E}[u(t_{k})]) \Delta t
+ \sigma B_k, \notag \\ & B_k:=B_{t_{k+1}}-B_{t_k} \sim \mathcal{N}(0, \Delta t),~ k \in \{0,\dots,N_T-1\}, \notag
\end{align}
and the associated cost function
\begin{align} \label{discrtized cost}
J(u) &= \mathbb{E} [\sum_{k=0}^{N_{T}-1} \ell(x(t_{k}),\mathbb{E}[x(t_{k})],u(t_{k}),\mathbb{E}[u(t_{k})]) \Delta t \notag\\ &+ \psi(x(t_{N_T}),\mathbb{E}[x(t_{N_T})])], 
\end{align}
In order to deal computationally with \eqref{discrtized cost},
we replace the expectations by empirical averages over $N$-dimensional sample of state trajectories, with initial conditions independently drawn from some distribution $\mu_0$. The same is done for the corresponding control trajectories, i.e., the cost functional becomes
\begin{align*}
&J^N(u) = \frac{1}{N}\sum_{i=1}^{{N}}\bigg(\psi(x^i({t_{N_T}}),\mu_x({t_{N_T}}))\\ 
&+\sum_{k=0}^{{N_T}-1} \ell(t_k,x^i({t_k}),\mu_x({t_k}),u^i(t_k), \mu_u({t_k})) \Delta t\bigg) ,\notag 
\end{align*}
where,
\begin{align*}
\mu_x({t_k)} &= \frac{1}{N}\sum_{j=1}^{{N}} x^j({t_k}),&
\mu_u({t_k)} &= \frac{1}{N}\sum_{j=1}^{{N}} u^j({t_k}). 
\end{align*}
The corresponding dynamics is given by
\begin{align}
\label{sde_dis_N} \notag
x^i({t_{k+1}})&= x^i({t_{k}}) + f(t_k,x^i({t_{k}}),\mu({t_{k}}),u^i(t_k), \mu_u({t_k}))) \Delta t \\
&+ \sigma B^i_k,~ x_0^i \sim \mu_0, \ B_k^i \sim \mathcal{N}(0, \Delta t), \\ \notag
k& \in \{0, \dots,N_T-1 \}, \ \ i \in \{ 1, \dots, N \}, \notag
\end{align}
and the problem we aim to solve is the following
\begin{align}
\label{eq:problem_NN_control_2}
\mathcal{P}_{\mathrm{CD}}:=
\left\{
\begin{array}{l}
\displaystyle \min_{u}~ J^N(u),\\
\text{subject to}~\eqref{sde_dis_N},~\text{and}~\\
x_0^i \sim 
\mu_0,~ i = 1,\dots, N. 
\end{array}
\right.
\end{align}

Therefore, solving the problem yields $N$ optimal controls and the corresponding optimal system states following the dynamics \eqref{eq:sde_dis_control},
\begin{align*}
\bar{u}^i &:=  (\bar{u}^i(t_{0},\cdot),\dots,\bar{u}^i(t_{N{_T}-1},\cdot) ),~i = 1,\dots, N,\\
%
\bar{x}^i &:= (\bar{x}^i(t_{0}),\dots,\bar{x}^i(t_{N{_T}-1}) ),~i = 1,\dots, N.
\end{align*}
For each time point we define the averages
\begin{align*}
\mu_{\bar{x}}(t_k)&= \frac{1}{N}\sum_{i=1}^{N} \bar{x}^i(t_k),~k = 1,\dots, N_T-1,\\
\mu_{\bar{u}}(t_k) &= \frac{1}{N}\sum_{i=1}^{N} \bar{u}^i(t_{k},\cdot),~k = 1,\dots, N_T-1.
\end{align*}

%
Now, we want to design a data-driven offline controller using all the previously generated optimal trajectories as training inputs for the neural network. We define the following loss function
\begin{align}
\label{eq:loss_function}
L^N(\theta) = \sum_{i = 1}^{N} \sum_{k=0}^{N{_T}-1} \left\| g_\theta(\cdot)-
     \bar{u}^i(t_k)
 \right\|,
\end{align}
where $\|\cdot\|$ denotes the Eucledian norm. The goal is to minimize $L^N(\theta)$ by searching for a suitable function (neural network) $g_{\bar{\theta}}(z) \in \mathcal{U}_{\mathbf{N}\mathbf{N}}$ parameterized by
$\bar{\theta} \in \Theta$. 
Solving this problem can be done efficiently with the
Stochastic Gradient Descent (SGD) algorithm, either by using TensorFlow or Pytorch in Python, for instance. Thus, the neural network training is achieved by solving the following optimization problem
\begin{align}
\label{eq:problem_NN_control}
\mathcal{P}_{\mathrm{trainNN}}:=
\left\{
\begin{array}{l}
\displaystyle \min_{\theta}~ L^N(\theta),\\
\text{subject to a compatible architecture}\\
\text{for the neural network}~ g_\theta,~\text{and}~\\
\bar{u}^i,~i = 1,\dots, N,~\text{given}.
\end{array}
\right.
\end{align}
The optimal weight and bias parameters are obtained from the solution of problem $\mathcal{P}_{\mathrm{trainNN}}$, i.e., $$ \bar{\theta} \in \arg\min_{\theta} L^N(\theta).$$

Next, we study the stability of the closed-loop neural network system composed of the stochastic system dynamics in \eqref{eq:sde_dis_control} and the trained neural network $g_{\bar{\theta}}$.

\section{Attacks and Stability for the Mean-Field-Type Neural Network}
\label{sec:attacks}


Once the neural mean-field-type control is trained, a natural question that arises is how to determine its performance. It has been reported, mainly, in classification problems (see \cite{Attacks_NN}) that a subtle modification in the input data might fool a well-trained neural network. These perturbed images are known as adversarial inputs or attacks. Here, we focus on the stabilization of the underlying stochastic system by the neural network using the concept of adversarial inputs. 
We make use of this concept in the analysis of the stability of the neural network mean-field-type control.
As we will see in Section \ref{sec:stability_num_example}, the neural network ensures stability properties within a small region $\mathcal{B}_\delta$. 
The adversarial inputs, denoted by $\check{x}(t_0) \notin \mathcal{B}_\delta$, are interpreted as attacks provided that the neural network cannot stabilize the stochastic system for such initial conditions. 

\subsection{Stochastic stability}
\label{sec:stochastic_stability}
%
%
Given $\gamma>0$, denote by $\mathcal{B}_{\gamma}$ the ball in $\mathbb{R}^d$ with center $0$ and radius $\gamma$: $\mathcal{B}_{\gamma} = \{ x : \|x\| \leq \gamma\}$.

The solution of the stochastic difference equation in \eqref{eq:sde_dis_control} is said to be stochastically stable if for every $\varepsilon \in (0,1)$ and every $r>0$ there exists $\delta:=\delta(\varepsilon,r)>0$ 
%
%
such that
$$\mathbf{P}( x({t_{k}}) \in \mathcal{B}_r~\text{for all}~k\in\{1,\dots,N_T-1\} ) \geq 1 - \varepsilon,$$ 
whenever $x(t_0) \in \mathcal{B}_\delta$. Otherwise, it is said to be stochastically unstable. This concept is general in the sense that we are not restricted to systems that are Markovian, stationary or ergodic.
\begin{figure}
    \centering
    \includegraphics[scale=1]{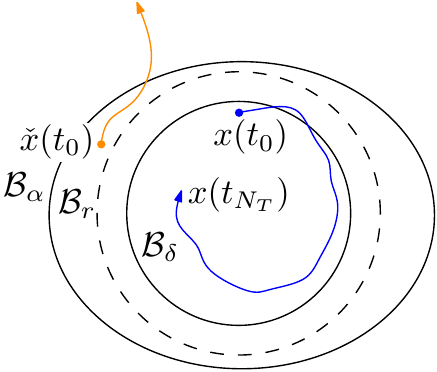}
    \caption{Illustration of two example trajectories, one starting within $\mathcal{B}_\delta$ (blue line) while the other starts with an  initial condition $\check{x}(t_0)$ within the set $\mathcal{B}_\alpha\setminus \mathcal{B}_\delta$ (orange line). Since the system is unstable then $\check{x}(t_0)$ is an adversarial attack}
    \label{fig:B_sets}
    \vspace{-0.6cm}
\end{figure}

For a given parameter $r$ defining the set $\mathcal{B}_r$, let $\mathcal{B}_\delta \subset \mathcal{B}_r$ be a set of initial conditions in which the closed-loop neural network control system is stochastically stable. In particular, if $\varepsilon$ is very close to zero, the trajectories starting in the set $\mathcal{B}_\delta$ remain inside $\mathcal{B}_r$ with probability very close to one.
Therefore, it is most likely not possible to find an adversarial initial condition in $\mathcal{B}_\delta$. This means, the probability to find an initial condition in $\mathcal{B}_\delta$ such that the trajectory leaves the set $\mathcal{B}_r$ and diverge is close to zero.

 Let us define the set $\mathcal{B}_\alpha$ where $\mathcal{B}_\delta \subset \mathcal{B}_\alpha$ and such that there exists an adversarial $\check{x}(t_0) \in \mathcal{B}_\alpha \setminus \mathcal{B}_\delta.$ All the previously defined sets are illustrated in Figure \ref{fig:B_sets}. 
 The numerical example presented in Section \ref{example}, will show the stochastic stability for different selection of the parameters $\varepsilon$ and $r$, and different sets characterizing the stability.

\subsection{Finding Adversarial Initial System State}

Let us fix a trained neural network $g_{\bar{\theta}}$. We are interested in finding an adversarial initial state in $\mathcal{B}_\alpha$ such that the following closed-loop trained neural network control system
\begin{align}
\label{eq:sde_dis_nn} \notag
x(t_{k+1}) &= x(t_{k}) + f(x(t_{k}),\mu_{x}(t_{k}), g_{\bar{\theta}}(t_k,\cdot), \mu_{g_{\bar{\theta}}}(t_{k}) ) \ \Delta t \notag
\\ &+ \sigma B_k,\ \ k \in \{0,\dots,N-1\},
\end{align}
is unstable. 
%
Such an adversarial initial state can be found by maximizing the cost functional that was used for training purposes, given the optimal neural network control $g_{\bar{\theta}}(\cdot)$, i.e.,
%
\begin{align*}
&J^N_{\bar{\theta}}\mathrm(x) = \frac{1}{N}\sum_{i=1}^{{N}}(\psi(x^i({t_{N_T}}),\mu_{x}({t_{N_T}}))\\ 
&+\sum_{k=0}^{{N_T}-1} \ell(t_k,x^i({t_k}),\mu_{x}({t_k}),g_{\bar{\theta}}^i(\cdot), \mu_{g_{\bar{\theta}}}({t_k})) \Delta t) ,\notag \\
&i \in \{0,\dots,N\}, \ k \in \{0,\dots,N_{T}-1\}.
\end{align*}
The goal is to solve the following optimization problem, 
\begin{align}
\label{eq:ad_problem_1}
\mathcal{P}_{\mathrm{AD}}:=
\left\{
\begin{array}{l}
\displaystyle \max_{x(t_0) \in \mathcal{B}_\alpha}~ J^N_{\bar{\theta}}(x(t_0)),\\
\text{subject to}~\eqref{eq:sde_dis_nn},\\
g_{\bar{\theta}}~\text{given}. 
\end{array}
\right.
\end{align}
%
The solution for  Problem \eqref{eq:ad_problem_1} defines an initial condition for the system state, $x(t_0) = \check{x}(t_0) \in \mathcal{B}_\alpha,$ such that the trained neural network $g_{\bar{\theta}}$ is unable to stabilize the system in \eqref{eq:sde_dis_control}. Problem \eqref{eq:ad_problem_1} can be solved by using, for example, the Projected Gradient Descent (PGD) algorithm as follows:
\begin{subequations}
\label{eq:SGD}
\begin{align}
{y}_{m+1} &= \mathrm{Proj}_{\mathcal{B}_\alpha} \left( {y}_{m} + \beta \cdot \nabla J^N_{\bar{\theta}}({y}_m) \right),\\
{y}_{0} &\in \mathcal{B}_\delta,~~m=0,\ 1,\dots
\end{align}
\end{subequations}
where $\mathrm{Proj}_{\mathcal{B}_\alpha}(\cdot)$ denotes the projection onto the set $\mathcal{B}_\alpha$, and $\beta$ denotes the step size for the gradient algorithm. 

The stopping condition is when the system state, following \eqref{eq:sde_dis_nn} with initial condition $x(t_0)={y}_m$, diverges, i.e., $|x(t_{N_T})| \to \infty$ (finite time escape that makes the state infinite). Thus, the adversarial is $\check{x}(t_0)=y_m$. 

\begin{figure}
    \centering
    \includegraphics[scale=1]{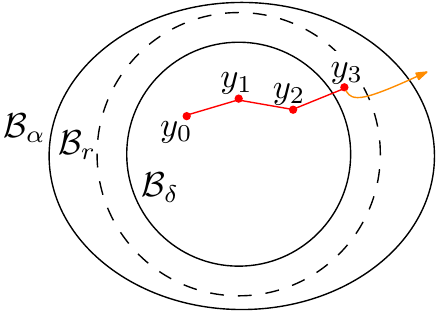}
    \caption{Evolution example of the PGD, which stops at $m=3$ given that the dynamics in \eqref{eq:sde_dis_nn} diverges with initial condition $x(t_0)=y_3$. Therefore, this is an adversarial $\check{x}(t_0)=y_3$. In the contrary, the dynamics in \eqref{eq:sde_dis_nn} are almost surely stochastically stable for initial conditions $x(t_0)\in \{y_0,y_1,y_2\}$.}
    \label{fig:evolution_SGD}
    \vspace{-0.6cm}
\end{figure}

Figure \ref{fig:evolution_SGD} shows an example for the algorithm  \eqref{eq:SGD}. The SGD finds an adversarial initial condition $y_m$ at step $m=3$ without requiring projection ($(y_{2} + \beta \cdot \nabla J^N_{\bar{\theta}}(y_2))$ does not leave the set $\mathcal{B}_\alpha$).
%
%
%

\subsection{Stability Improvement Using Adversarial Data}
\label{sec:improvingNN}

The neural network $g_{\bar{\theta}}$ in Section \ref{sec:NN} is trained by using initial conditions $x_0^j$, $j = 1,\dots, N$. The stability for the closed-loop neural network system is characterized by finding a number of adversarial initial states solving Problem \eqref{eq:ad_problem_1}, as explained in section VI. Now, the goal is to improve the stability of the neural network controller by enlarging the original training set to include, say, $\check{N} \in \mathbb{N}$ optimal state trajectories that are obtained from adversarial initial conditions, i.e., we generate adversarial inputs $(\check{x}^q(t_0) \in \mathcal{B}_\alpha,\ q=1,\ldots,\check{N})$ by solving \eqref{eq:ad_problem_1}, then we find the corresponding optimal controls by solving Problem \eqref{eq:problem_NN_control_2}, and we denote them by
%
\begin{align*}
\bar{v}^q &:=  (\bar{v}^q(t_{0},\cdot),\dots,\bar{v}^q(t_{N{_T}-1},\cdot) ),~q = 1,\dots, \check{N},
\end{align*}
The next step is to use this augmented data to retrain the neural network $g_\theta(\cdot)$. Thus, we aim to minimize the following loss function
\begin{align}
\label{eq:loss_function_improved}
L&^{N}_2(\theta)  \\
&=  \sum_{k=0}^{N{_T}-1} \Bigg[ \sum_{j = 1}^{N} \left\| g_\theta(\cdot)-
     \bar{u}^j(t_k)  
 \right\| + \sum_{q = 1}^{\check{N}} \left\| g_\theta(\cdot)-
     \bar{v}^q(t_k)  
 \right\|
 \Bigg]. \notag
 \end{align}
The incorporation of this larger set of data should improve the approximation of the optimal control input without requiring to modify the architecture of the neural network.

In the coming section, we present a numerical example consisting of a control design for an unstable dynamical stochastic system. 

\section{Numerical Example}\label{example}

For illustrative purposes,  we train two neural network mean-field-type controllers with quadratic cost functional and linear dynamics as follows:
\begin{align*}
\ell(t_k&,x({t_k}),\bar{x}({t_k}),u_\theta(\cdot)) \Delta t \\
&= q_1({t_k}) x({t_k})^2 + r_1({t_k}) u_\theta(\cdot)^2\\
&+ q_2({t_k}) (x({t_k})-\mathbb{E}[x({t_k})])^2 + r_2({t_k}) (u_\theta(\cdot)-\mathbb{E}[u_\theta(\cdot)])^2,\\
\psi(&x({t_{N_T}}),\bar{x}({t_{N_T}})) \\
&= q_1({t_{N_T}}) x({t_{N_T}})^2 + q_2({t_{N_T}}) (x({t_{N_T}})-\mathbb{E}[x({t_{N_T}})])^2,
\end{align*}
and
\begin{align}
\label{eq:example_dynamics}
    x({t_{k+1}}) &= a_1 x({t_{k}}) + a_2 \mathbb{E}[{x}({t_{k}})] \notag\\
    &+ b_1 u_\theta(\cdot) + b_2 \mathbb{E}[{u}_\theta(\cdot)] + \sigma B_k,
\end{align}
%
where $a_1 = 2$, $a_2= 1$, $b_1= 1$, $b_2 = 2$, and $\sigma=1$. The parameters of the cost functional are: $q_1 = 20$, $q_2 = 10$, $r_1 = 200$, $r_2 = 100$. The time horizon is fixed to be $N_T = 15$ with $\Delta t = 1/20$. For this linear example, it is not needed to compute the expectation of the system states in the empirical form as an average of multiple trajectories. Instead, it is possible to compute the evolution of the expected system state $\mathbb{E}[x({t_{k}})]$.\\
%
The system dynamics in \eqref{eq:example_dynamics} is \textit{unstable}, this can be seen from the evolution of the expected state $\mathbb{E}[x({t_{k+1}})]$ given that $(a_1 + a_2)>1$.
%
\subsection{Training Stage}
%


\begin{figure}[t!]
	\begin{center}
		\resizebox{\columnwidth}{!}{
			\begin{tabular}{cc}
				&\textbf{Neural Network 1} \\ 
				\rotatebox{90}{\hspace{1.5cm} $\mathbb{E}[x(t_k)]$} & \includegraphics[scale=0.3]{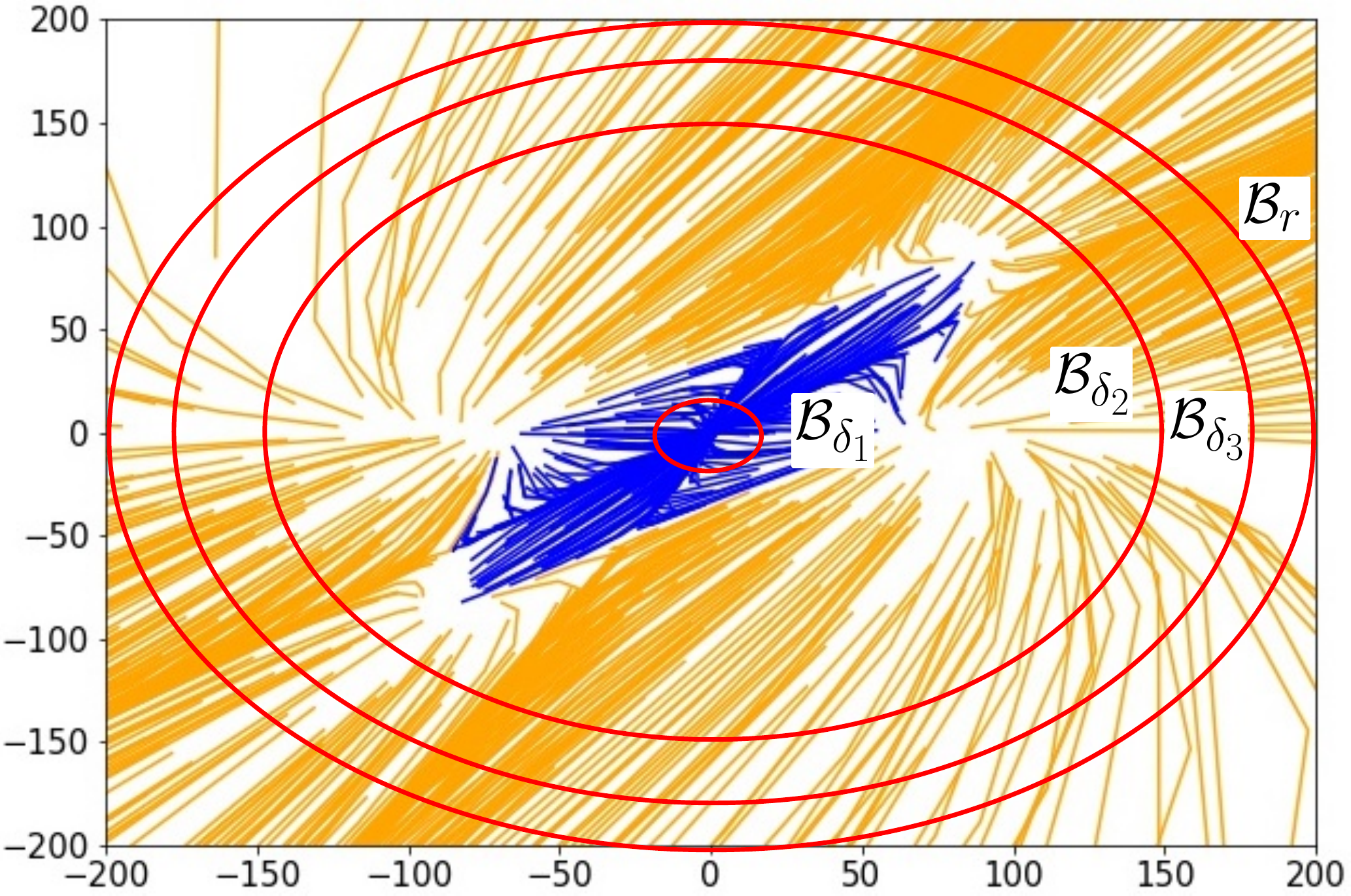}\\
				& ~~$x(t_k)$ \\ 
				& \textbf{Neural Network 2} \\ 
				\rotatebox{90}{\hspace{1.5cm} $\mathbb{E}[x(t_k)]$} & \includegraphics[scale=0.3]{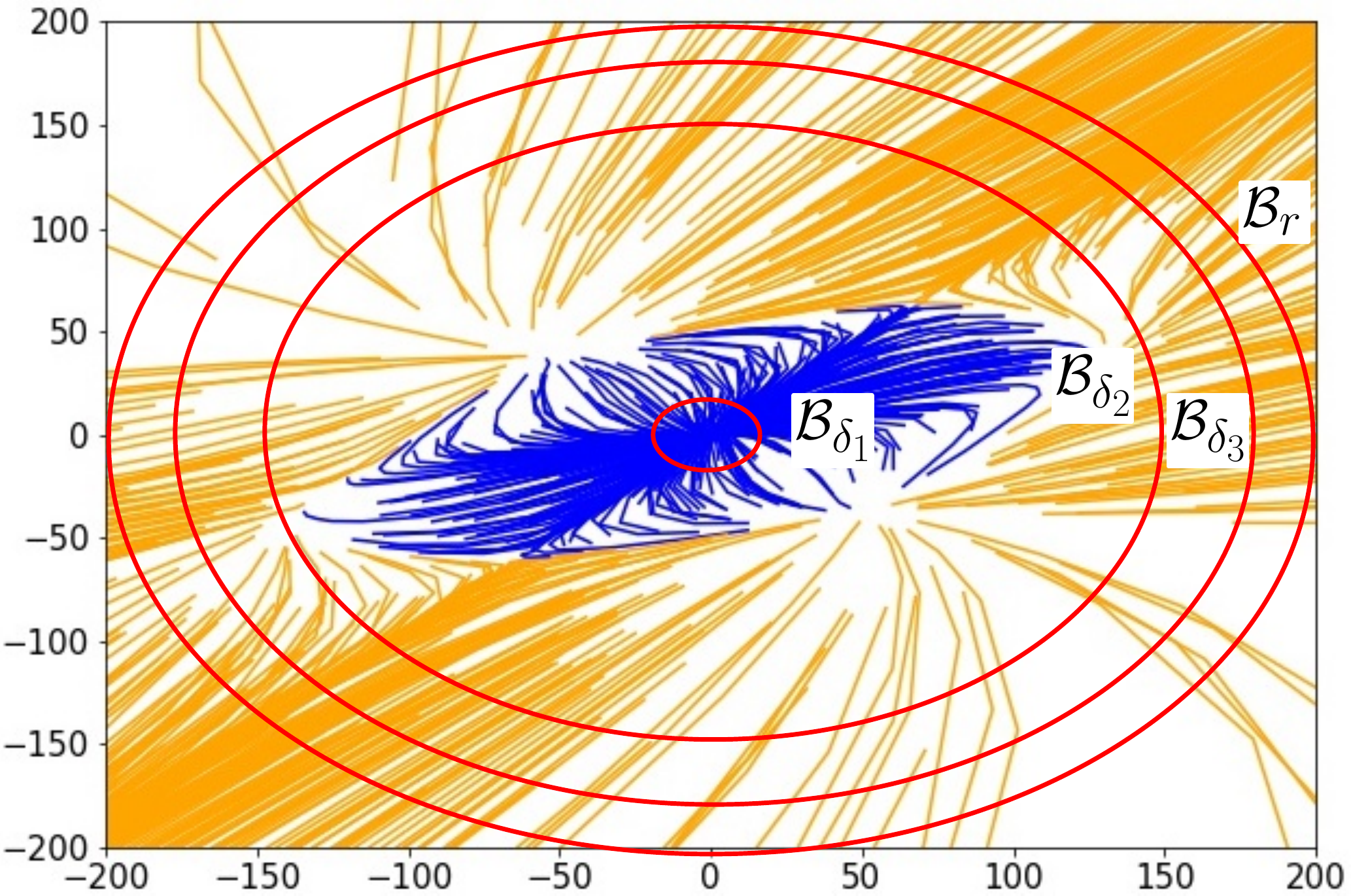} \\
				& ~~$x(t_k)$ \\ 
				& \textbf{Improved NN 1 (Adversarial training)}\\  
				 \rotatebox{90}{\hspace{1.5cm} $\mathbb{E}[x(t_k)]$} & \includegraphics[scale=0.3]{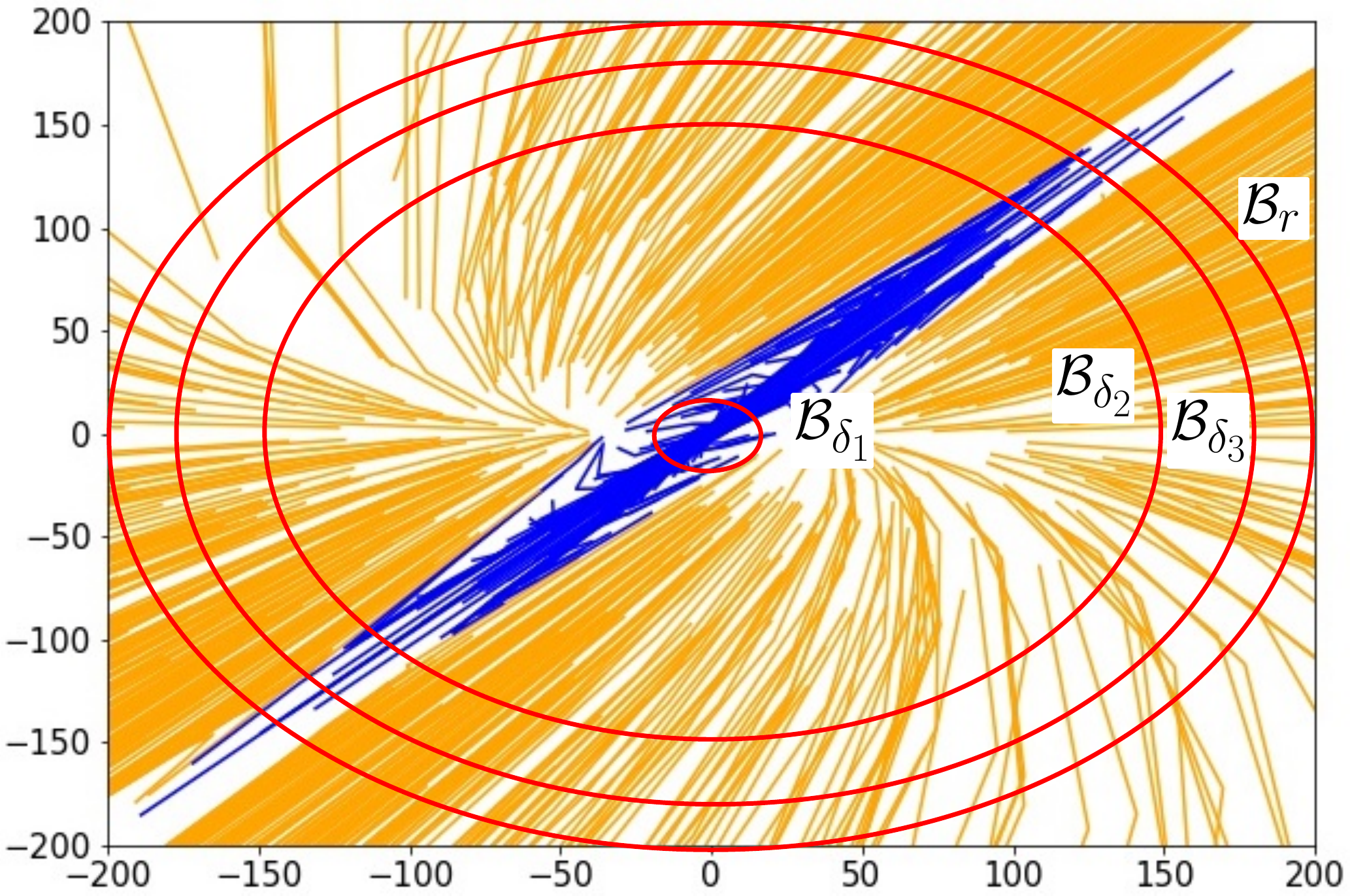}\\				
				& ~~$x(t_k)$\\
			\end{tabular}
		}
	\end{center}
	\vspace{-0.3cm}
	\caption{Stabilization of the stochastic system by the trained neural network mean-field-type control. Blue trajectories are stable, and orange trajectories are unstable.}
	\label{fig:comparison}
	\vspace{-0.5cm}
\end{figure}

We consider two architectures for two different neural networks whose output  mappings are 
$g^1_{\theta}, g^2_{\theta}:\mathbb{R}^3 \to \mathbb{R}^2$, respectively, and with the same input given by
$z({t_k}) = [x({t_k}) \quad \mathbb{E}[x({t_k})] \quad B_k]^\top$ 
and output $g^1_\theta(\cdot), g^2_\theta(\cdot) := [u_\theta(\cdot) \quad \mathbb{E}[{u}_\theta(\cdot)]]^\top.$ The characteristics of the neural networks are shown in Table \ref{tab:NN_architectures}.
\begin{table}[t!]
    \centering
    \caption{Architecture for the two different tested neural networks}
    \label{tab:NN_architectures}
    \resizebox{\columnwidth}{!}{
    \begin{tabular}{l|cc}
    \hline
    \textbf{Parameter} & \textbf{Neural Network 1} & \textbf{Neural Network 2} \\
    \hline
    \hline
    \textbf{Layers} & $3$ & $5$ \\
    \textbf{Total  Neurons} & $6$ & $106$ \\
    \textbf{Neurons per Layer} & $\{2,2,2\}$ & $\{2,2,50,50,2\}$ \\
    \textbf{Activation Functions} & $\{\text{lin},\text{tanh},\text{lin}\}$ & $\{\text{lin},\text{tanh},\text{tanh},\text{tanh},\text{lin}\}$\\
    \hline
    \end{tabular}
    }
    \vspace{-0.5cm}
\end{table}
The layers are characterized by the following parameters for a given architecture $g^1_\theta$ or $g^2_\theta$:
\begin{align*}
 \theta|g^1_\theta &:= \{W^{(0)},b^{(0)},W^{(1)},b^{(1)},W^{(2)},b^{(2)} \},\\
 \theta|g^2_\theta &:= \{W^{(0)},b^{(0)},W^{(1)},b^{(1)},W^{(2)},b^{(2)},\dots,W^{(4)},b^{(4)} \},
\end{align*}
and the output of the neural networks is as follows:
\begin{align*}
g^1_{\theta}&=\phi^{2} \circ \phi^{1}  \circ \phi^0(z),&
g^2_{\theta}&=\phi^{4} \circ \phi^{3} \circ \phi^{2} \circ \phi^{1}  \circ \phi^0(z).
\end{align*}
Next section discusses the stability properties of the neural network closed-loop dynamics. 

\subsection{Stability Assessment}
\label{sec:stability_num_example}

We are now interested in evaluating the stochastic stability for the neural network closed-loop dynamics
\begin{align}
\label{eq:nn_closedloop}
    x({t_{k+1}}) &= a_1 x({t_{k}}) + a_2 \mathbb{E}[{x}({t_{k}})] + 
    \sigma B_k\notag \\
    &+ 
    [b_1 \quad b_2] \cdot g_{\bar{\theta}}(x({t_{k}}),\mathbb{E}[x({t_{k}})],B_k),
\end{align}
and for the two considered neural networks, i.e., $g_{\bar{\theta}}(\cdot) := \{g^1_{\bar{\theta}}(\cdot), g^2_{\bar{\theta}}(\cdot)\}$.
Note that the stability analysis for the stochastic system dynamics in \eqref{eq:nn_closedloop} is involved given that the functions 
\begin{align*}
g^1_{\theta}(\cdot) &= \phi^{2} \circ \phi^{1}  \circ \phi^0([x(t_k) \quad \mathbb{E}[x(t_k)] \quad B_k]),\\
g^2_{\theta}(\cdot) &= \phi^{4} \circ\phi^{3} \circ \phi^{2} \circ \phi^{1}  \circ \phi^0([x(t_k) \quad \mathbb{E}[x(t_k)] \quad B_k]),
\end{align*}
are non-linear according to the selected activation functions.  

To this end, we will computationally characterize the set $\mathcal{B}_{\delta}$ for the stability using the first neural network (see Table \ref{tab:NN_architectures}) and for the given values $r$ and $\varepsilon$. Such sets are determined by
empirically computing the following probability:
\begin{align}
\label{eq:probability}
    \mathbf{P}( x(t_{k})  \in  \mathcal{B}_r, \forall~k  \in \{1,...,N_T-1\}, x(t_0) \in  \mathcal{B}_\delta ) .
\end{align}
Let us consider the following three scenarios 
\begin{itemize}
	\item \textbf{Scenario 1:} $r=200$, and $\varepsilon_1 \approx 0$,
	\item \textbf{Scenario 2:} $r=200$, and $\varepsilon_2=0.55$,
	\item \textbf{Scenario 3:} $r=200$, and $\varepsilon_3=0.7$,
\end{itemize}
which define the sets $\mathcal{B}_{\delta_1}(r,\varepsilon_1)$, $\mathcal{B}_{\delta_2}(r,\varepsilon_2)$, and $\mathcal{B}_{\delta_3}(r,\varepsilon_3)$, respectively. We test $1000$ trajectories corresponding to random initial states in order to find the values 
\begin{align*}
\delta_1&=20,& 
\delta_2&=150,&
\delta_3&=180. 
\end{align*}
Table \ref{tab:comparison} shows details in its first column corresponding to the first neural network and the established values for $r$ and $\varepsilon$.
Figure \ref{fig:comparison}(a) shows the evolution of the system state $\bar{x}(t_k)$ and its expectation $\mathbb{E}[\bar{x}(t_k)]$ according to \eqref{eq:example_dynamics} and using the optimal control input ${u}_{\bar{\theta}}(t_k)$ computed by means of the neural networks  $g^1_\theta$. 
\subsection{Stability Comparison between Architectures}
In order to compare the stability of the neural network closed-loop dynamics for the two different architectures (see Table \ref{tab:NN_architectures}), we compute the probability in \eqref{eq:probability} for the two neural network architectures using the same sets $\mathcal{B}_r$, $\mathcal{B}_{\delta_1}(r,\varepsilon_1)$, $\mathcal{B}_{\delta_2}(r,\varepsilon_2)$, and $\mathcal{B}_{\delta_3}(r,\varepsilon_3)$. Table \ref{tab:comparison} shows the stability comparison between the neural networks 1 and 2, showing better stability properties when using the second neural network that is composed of more layers and total neurons. 
\begin{table}[t!]
	\caption{Neural network comparison according to their stability.}
	\label{tab:comparison}
	\begin{center}
		\begin{tabular}{c|ccc}
			\hline
			& \multicolumn{3}{c}{$\mathbf{P}( x(t_{k}) \in \mathcal{B}_r\text{, for all}~k,\ x(t_0) \in \mathcal{B}_\delta )$} \\
			\textbf{Scenario} $(r=200)$	& \textbf{NN 1} & \textbf{NN 2} & \textbf{Improved NN 1} \\
			\hline
$\mathcal{B}_{\delta_1}$, $\delta_1=20$ & 1 & 1 & 1 \\
$\mathcal{B}_{\delta_2}$, $\delta_2=150$ & 0.45 & 0.557 & 0.464 \\
$\mathcal{B}_{\delta_3}$, $\delta_3=180$ & 0.3 & 0.449 & 0.354 \\
			\hline
		\end{tabular}
	\end{center}
\end{table}

\subsection{Stability Improvement}

We have observed in Section \ref{sec:stability_num_example} that the modification of the neural network architecture leads to an improvement in the stochastic stability of the closed-loop neural network. For example, the performance of the neural network with $106$ neurons and $5$ layers exhibited a higher stability probabilities than the neural network with just $6$ neurons and $3$ layers. 

An alternative to improve the performance of the closed-loop neural network consists of enlarging the training set to re-train the neural network using adversarial inputs as presented in Section \ref{sec:improvingNN}. We improve the neural network with $3$ layers and $6$ neurons following this methodology (neural network $1$ in Figure \ref{fig:comparison}). To this end, we generate $\check{N}=500$ adversarial inputs and re-train the neural network $1$ minimizing the loss function in \eqref{eq:loss_function_improved}.

Figure \ref{fig:comparison}(c) shows an improvement of the stability properties for low-variance initial conditions. For instance, when $x_0=\mathbb{E}[x_0]$, the stability of the neural network closed-loop using the neural network $1$ shows convergence for values in the range $-80<x_0<80$, whereas for the improved neural network $1$, this range is enlarged to be $-190<x_0<190$. In addition, we compute the probability \eqref{eq:probability} for the closed-loop dynamics using the improved neural network $1$ as shown in Table \ref{tab:comparison}. It can be seen an improvement with respect to the neural network $1$.

\color{black}
\section{Concluding Remarks and Future Directions}

We have presented a data-driven mean-field-type control via neural networks. We have studied the stability of the closed-loop neural network control system. This is done by using, first, a simple two-stage method to train the neural network, then characterizing the basin of attraction by means of adversarial inputs which also, in a sense, validate the feasibility of the solutions obtained from the neural network. Furthermore, we proposed a way to improve the robustness and the stability of the approximated solutions by adversarial training. Finally, we numerically compared two different neural-network architectures. The results suggested that more complex neural networks might lead (but not guaranteed) to more robust mean-field-type control, i.e., a bigger invariant forward basin of attraction has been observed with deep learning (with more hidden layers in the architecture). Moreover, we showed that an adversarial training can significantly enlarge the basin of attraction and thus the stability and the robustness for less complex (smaller) neural-network architectures (without the need of modifying the architecture).

As further work, we propose to extend the results presented in this paper to the game theoretic case, i.e., study stability of neural networks mean-field-type games using adversarial attacks. In addition, to address several game solution concepts such as non-cooperative, zero-sum, Stackelberg, hierarchical, and Berge games, among others. Moreover, the rigorous mathematical characterization of the stability sets for neural networks
solving mean-field type control problem, which were computationally estimated in this work, is an open theoretical problem
to be considered in the future.

\bibliographystyle{unsrt}
\bibliography{references}

\end{document}